\documentclass[reqno]{amsart}
\usepackage{epstopdf}
\usepackage{subfigure}

\usepackage{natbib}
\bibpunct[, ]{(}{)}{;}{a}{}{,}
\usepackage{hyperref}
\usepackage{bbm, comment}
\usepackage{color}

\theoremstyle{plain}
\newtheorem{theorem}{Theorem}[section]
\newtheorem{lemma}[theorem]{Lemma}
\newtheorem{corollary}[theorem]{Corollary}
\newtheorem{proposition}[theorem]{Proposition}

\theoremstyle{definition}

\theoremstyle{remark}

 \DeclareMathAlphabet{\mathpzc}{OT1}{pzc}{m}{it}
 \DeclareMathAlphabet{\mathsfsl}{OT1}{cmss}{m}{sl}

  \newcommand{\FH}{\mathfrak{H}}

\newcommand{\dif}{\mathrm{d}}

\newcommand{\R}{\mathbb{R}}

\newcommand{\abs}[1]{\left\vert#1\right\vert}
\newcommand{\set}[1]{\left\{#1\right\}}

\newcommand{\norm}[1]{\left\Vert#1\right\Vert}
 \newcommand{\innp}[1]{\langle {#1}\rangle}

\newcommand{\E}{\mathbb{E}}

\newcommand\sgn{\mathrm{sgn}}

\newcommand{\blue}{\textcolor[rgb]{0.00,0.00,1.00}}
\newcommand{\red}{\textcolor[rgb]{1.00,0.00,0.00}}

\begin{document}
\title[asymptotic expansion of the norm of $e^{-\abs{t-s}}\mathbf{1}_{\set{0\le s,t\le T}}$]{An asymptotic expansion of the norm of $e^{-\abs{t-s}}\mathbf{1}_{\set{0\le s,t\le T}}$ in the canonical Hilbert space of fractional Brownian motion}
\author[Y. Chen]{Yong Chen}
\address{School of Big Data, Baoshan University, Baoshan, 678000, Yunnan, China.}
\email{zhishi@pku.org.cn}
\begin{abstract}
Using the inner product formula of the canonical Hilbert space of fractional Brownian motion on an interval $[0,T]$ with Hurst parameter $H\in (0,1)$ given by Alazemi et al., we show the asymptotic expansion of the norm of $f_T(s,t):=e^{-\abs{t-s}}\mathbf{1}_{\set{0\le s,t\le T}}$  up to the term $T^{4H-4}$. As applications, we show that the existence of the oblique asymptote of  the norm $\frac12\|f_T\|^2_{{\FH}^{\otimes2}}$ if and only if $H\in (0,\frac12]$ and that we obtain a sharp upper bound of the difference $\abs{\frac{1}{2 {T}} \norm{f_T}_{\FH^{\otimes 2}}^2-\sigma^2}$ for $H\in (0,\frac34)$ which implies two significant estimates concerning to an ergodic fractional Ornstein-Uhlenbeck  process, where $\sigma^2$ is the slope of the oblique asymptote.\\
{\bf Keywords :} fractional Brownian motion; canonical Hilbert space; fractioanl Ornstein-Uhlenbeck (OU) process; second chaos random variable.\\
{\bf MSC 2020:} 60G22; 41A60; 60H07.

\end{abstract}
\maketitle

\section{Introduction}
The inner product formula of the canonical Hilbert space $\FH$ of fractional Brownian motion (fBm) on an interval $[0,T]$ with Hurst parameter $H\in (0,1)$ is well known. Especially, when $H\in (0,\frac12)$, we have 
	\begin{align}
		 \langle\phi,\psi\rangle_{\FH}=\langle K^*_H\phi(t),K^*_H\psi(t)\rangle_{L^2([0,T])}, \label{suanzi bianhuanfa}
	\end{align}where 
operator~$K^*_H$
	\begin{align*}
		 (K^*_H\phi)(t)=K_H(T,t)\phi(T)+\int^T_t\frac{\partial K_H}{\partial s}(s,t)[\phi(s)-\phi(t)]\dif s,
	\end{align*} and 
    	\begin{align*}
		K_H(t,s)=c_H\Big[\Big(\frac ts\Big)^{H-\frac12}(t-s)^{H-\frac12}-(H-\frac12)s^{\frac12-H}\int^t_su^{H-\frac32}(u-s)^{H-\frac12}\dif u\Big],
	\end{align*}
see~[12] 
for example. It is clear that \eqref{suanzi bianhuanfa} is not easy to be used to compute the inner product.

In [11], 
the following formula provides a computation for the inner product of two functions in the canonical Hilbert space $\mathfrak{H}$:
\begin{align}\label{zh-n-hu}
    \langle f,\,g \rangle_{\FH}&=-\int_{[0,T]^2}  f(t) g'(s)  \frac{\partial R(s,t)}{\partial t}  \dif t \dif s,
\end{align} 
where $$R(s,t)=\frac{1}{2}[s^{2H}+t^{2H}-\abs{s-t}^{2H}]$$ is the covariance function of the fBm, and the derivative $g'(s)$ can understood as the distribution derivative (see [3]). 
In 
[1], they show that this formula can be simplified using the fact that the term $t^{2H-1}$ in the partial derivative
\begin{align*}
    \frac{\partial R(s,t)}{\partial t}=H\big[t^{2H-1} -\abs{t-s}^{2H-1}\sgn(t-s)\big] 
\end{align*} does not contribute to the integration value given by \eqref{zh-n-hu}. This simplified formula is summarized in the following Proposition \ref{st1-thm-original}.

 Denote by $\mathcal{V}_{[0,T]}$ the set of functions of bounded variation on $[0,T]$, and by $\mathcal{B}([0,T ])$ the
Borel $\sigma$-algebra on $[0,T]$. 
For each $g\in \mathcal{V}_{[0,T]}$, $\nu_{g}$ is the restriction to $([0,T ], \mathcal{B}([0,T ]))$ of the Lebesgue-Stieltjes signed measure $\mu_{g^0}$ on $\left(\R,\mathcal{B}(\R)\right)$, where $g^0(x)$ is defined by 
\begin{equation*} 
  g^0(x)=\left\{
  \begin{array}{ll}
    g(x), & \quad \text{if } x\in [0,T] , \\
    0,    & \quad \text{otherwise}.
  \end{array}
  \right.
\end{equation*}
\begin{proposition}\label{st1-thm-original}
Let $H \in (0, \frac{1}{2})\cup(\frac12,1)$. For any two functions in the set $\mathcal{V}_{[0,T]}$, their inner product in the Hilbert space $\mathfrak{H}$ can be expressed as
 \begin{align} 
\langle f,\,g \rangle_{\FH}
&=H \int_{[0,T]^2}  f(t)  \abs{t-s}^{2H-1}\sgn(t-s) \dif t \nu_{g}(\dif s),\quad  \forall f,\, g\in \mathcal{V}_{[0,T]}. \label{innp fg3-0}
\end{align}    
If $  g'(\cdot) $ is interpreted as the distributional derivative of $g(\cdot)$, then the formula \eqref{innp fg3-0} admits the following representation:
 \begin{align*} 
\langle f,\,g \rangle_{\FH}
&=H \int_{[0,T]^2}  f(t) g'(s) \abs{t-s}^{2H-1}\sgn(t-s) \dif t  \dif s,\quad  \forall f,\, g\in \mathcal{V}_{[0,T]}.
\end{align*}  
\end{proposition}

Although \eqref{innp fg3-0} is restricted to functions in the set $\mathcal{V}_{[0,T]}$, Proposition~\ref{st1-thm-original} is still very useful because the set $\mathcal{V}_{[0,T]}$ is large enough for some applications, see 
[1]. In the present paper, we will follow the same idea to show the asymptotic expansion of the norm of the function of two variables $f_T$ (see \eqref{ft bsD}). As far as we know, the following asymptotic expansions, both \eqref{diyigedingli-1} and \eqref{34 case zkai}, are novel. Recall the notation $g(u)=O(f(u))$ means that there exist constants $M$ and $c$ such that the real-valued functions $ f,\,g$ satisfies $\abs{g(u)}\le M f(u)$ for all $u>c$.  
\begin{theorem}\label{qq2b}
Let~
 $H\in(0,1)$. For the function of two variables
	\begin{equation} \label{ft bsD}
	 f_T(t,s)=e^{- |t-s|}\mathbf{1}_{\{0\leq s,t\leq T\}}
	\end{equation}
 in the symmetric tensor space $\FH^{\odot 2}$, the equality 
 \begin{align}  		 \|f_T\|^2_{{\FH}^{\otimes2}}&=2{H^2}\Gamma^2(2H)\left[(4H-1)(1-\frac{1}{\cos(2H \pi)})(T-\frac{4H-1}{2})-(2H-1)(2H+1) \right]\notag\\
    &+\frac{4(2H-1)}{4H-3}T^{4H-2}-\frac{8(2H-1)^2}{4H-3}T^{4H-3}+O(T^{4H-4}) \label{diyigedingli-1}
 \end{align}	 
	holds when $H\neq \frac34$; and 
\begin{align}
 \|f_T\|^2_{{\FH}^{\otimes2}}&= \frac{9}{8}\left[(T-1)\log T+(c+\frac{\pi-3}{2})T+1-c-\frac{13\pi}{16}+\frac{1}{T}\right]+O(T^{-2})\label{34 case zkai}
\end{align}
    holds when $H=  \frac34$, where $c=2\log 2+\gamma$ and $\gamma$ is the Euler's constant. Moreover, the term $\frac{4H-1}{-\cos(2H\pi)}$ in \eqref{diyigedingli-1} is illustrated as its limit $\frac{2}{\pi}$ when $H=  \frac14$. 
\end{theorem}

    Denote 
\begin{align}
a&=H\Gamma(2H),\label{a bds}\\
   \sigma^2_H&=  (4H-1)(1-\frac{1}{\cos(2H \pi)}),\\
\sigma^2&=a^2\sigma_H^2= H^2\Gamma^2(2H) (4H-1)(1-\frac{1}{\cos(2H \pi)}).   \label{sigma2bds}
\end{align}
The following three corollaries can be obtained directly from Theorem~\ref{qq2b}.
 \begin{corollary}The following  limit exists if and only if $H\in (0,\frac34)$:
     \begin{eqnarray}\label{jixiansig2}
         \lim_{T\to\infty}\frac{1}{2 {T}} \norm{f_T}_{\FH^{\otimes 2}}^2=\sigma^2.
     \end{eqnarray}
 \end{corollary}When $H\in (\frac12,\frac34)$, \eqref{jixiansig2} has appeared in [10]. 
 When $H\in (0,\frac12)$, \eqref{jixiansig2} has appeared in [11] 
 by means of \eqref{zh-n-hu} where they do a very hard computation.
\begin{corollary}  
  There exists an oblique asymptote for  $\frac12\|f_T\|^2_{{\FH}^{\otimes2}}$ as a function of $T$ if and only if $H\in (0,\frac12]$.  Let $a$ and $\sigma^2$ be given by \eqref{a bds} and \eqref{sigma2bds}, respectively. The slope and the $y$-intercept of the oblique asymptote are equal to $\sigma^2$ and $-\frac{4H-1}{2}\sigma^2-(2H-1)(2H+1)a^2$, respectively.
\end{corollary}    
 \begin{corollary}\label{corjjx}
     When $H\neq \frac34$, there exists a positive constant $C$ independent of $T$ such that  
    \begin{align}\label{upper bound}
    \abs{\frac{1}{2 {T}} \norm{f_T}_{\FH^{\otimes 2}}^2-\sigma^2}\le \frac{C}{T^{1\wedge (3-4H)}}.
\end{align} 
 \end{corollary}
Although both the existence of the oblique asymptote of $\frac12\|f_T\|^2_{{\FH}^{\otimes2}}$ and the upper bound \eqref{upper bound} are obtained in 
\blue{[3]} when $H\in (0,\frac12)$,  the exact value of the slope of the oblique asymptote is borrowed from 
[11], and the exact value of the $y$-intercept of the oblique asymptote has not appeared in the literature before. This is due to that a very complicated inner product formula of the canonical Hilbert space of fBm is used  in 
\blue{[3]}.

The paper is organized as follows. In Section~\ref{sec.2}, we apply Theorem~\ref{qq2b} to an ergodic fractional Ornstein-Uhlenbeck (OU) process. In Section~\ref{sec.3}, we show Theorem~\ref{qq2b} using Proposition~\ref{st1-thm-original}.
 Several technical inequalities are provided in the Appendix.
 The symbol $C$ stands for a generic constant, whose value can change from one line to another.   
 
\section{Several applications to the ergodic fOU process}\label{sec.2}

In this Section, we will apply Theorem~\ref{qq2b} to an ergodic fractional OU process
 $\eta:=(\eta_t)_{t\in [0,T]}$,
which is defined as the solution to the stochastic differential equation
\begin{align}
   \dif \eta_t &=-\theta \eta_t \dif t +\dif  B^H_t,  \label{OU dingyi duibi}
\end{align} where $\theta$ is a positive constant and $B_t^H$ is the fBM in the interval $[0,T]$. The following Propositions~\ref{pp.2.1} and \ref{pp.2.2} and Corollary~\ref{coro.2.3} are significant to the problems of the parameter estimate of \eqref{OU dingyi duibi}, see 
[2-6; 8-9] for example.
  
{ Without loss of generality, we assume $\theta=1$ in the rest of this section.}
\subsection{an estimate of the difference between $\frac{\norm{f_T}_{\FH^{\otimes 2}}^2}{2 T}$  and $\sigma^2_H b_T^2$}

If we assume $\eta_0 = 0$, the solution to \eqref{OU dingyi duibi} has an integral representation 
 \begin{align} 
   \eta_t = \int_0^t e^{-\theta(t-u)}\dif B^H_u.
 \label{zeta T jifen}
\end{align} Denote  $h_t(u)=e^{u-t}\mathbf{1}_{[0,t]}(u)$. Then we have
\begin{align*}
    b_T:=\frac{1}{T}\int_0^T\, \E[\eta_t^2] \dif t=\frac{1}{T}\int_0^T\,  \norm{h_t}^2_{\FH}\dif t. 
\end{align*}

The following proposition is used to show an Berry-Esséen bound of an least squared estimate for the parameter $\theta$ of the fOU process $\eta$ in 
[4]. However, a very complicated inner product formula of fBm is used there. In the present paper, these results can be implied from Corollary~\ref{corjjx}. 
\begin{proposition}\label{pp.2.1}
When $H\in (0,\frac34)$, there exists a positive constant $C$ independent of $T$ such that 
    \begin{align*}
       \abs{\frac{\norm{f_T}_{\FH^{\otimes 2}}^2}{2 T} -\sigma^2_Hb^2_T} \le \frac{C}{{T^{1\wedge (3-4H)}}}.
    \end{align*}
\end{proposition}
\begin{proof}
 It follows from Proposition~\ref{st1-thm-original} that
    \begin{align*}
        b_T&=\frac{H}{T}\int_0^T \dif t\int_{[0,t]^2} e^{u-t+v-t}\left(\mathbf{1}_{[0,t](v)}+\delta_{0}(v)-\delta_t(v)\right) \abs{u-v}^{2H-1}\sgn(u-v)\dif u\dif v \\
        &=\frac{H}{T}\int_0^T \dif t\int_0^t e^{u-t}\left[e^{-t}u^{2H-1} +(t-u)^{2H-1}\right]\dif u, 
    \end{align*}
    which implies that there exists a positive constant $C$ independent of $T$ such that 
    \begin{align*}
        \abs{b_T-a}\le \frac{C}{T},
    \end{align*} see Lemma 3.7 of [3]. 
    Hence, we have that for $T$ large enough,
   \begin{align*}
        \abs{b_T^2-a^2}=\abs{(b_T-a)(b_T+a)}\le \frac{3aC}{T}.
    \end{align*} 
It follows from Corollary~\ref{corjjx} that  there exists a positive constant $C$ independent of $T$ such that
\begin{align*}
     \abs{\frac{\norm{f_T}_{\FH^{\otimes 2}}^2}{2 T} -\sigma^2_Hb^2_T} \le  \abs{\frac{\norm{f_T}_{\FH^{\otimes 2}}^2}{2 T} -\sigma^2 }+\sigma_H^2\abs{b^2_T-a^2}\le \frac{C}{T^{1\wedge (3-4H)}}.
\end{align*}
\end{proof}
\subsection{An upper bound the second moment of a second chaos random variable}
Let $\eta$ be given as in \eqref{zeta T jifen}. The following {second chaos random variable} with respect to the fBm $B^H$
    \begin{align} 
W_T=\frac{1}{\sqrt {T}}\int_{0}^T \big( \eta_t^2 -\E[\eta_t^2] \big)\dif t\label{wt dyi} 
\end{align}
plays a crucial role for the study of the Berry-Esséen bounds of both the least squared estimate and the moment estimate for the parameter $\theta$ of the fOU process $\eta$ in 
[5]. Especially, the following proposition is the start point there.
\begin{proposition}\label{pp.2.2}
When $H\in (0,\frac34)$, there exists a positive constant $C$ independent of $T$ such that the following estimates hold.
  \begin{align}\label{erjielieliang}
\abs{\E [W_T^2 ]-\sigma^2}&\le \frac{C}{{T}^{\gamma}},
\end{align}   where 
\begin{equation*}
\gamma=\left\{
      \begin{array}{ll}
 \frac12, & \quad \text{if }  H\in (\frac12,\,\frac58],\\
3-4H, &\quad \text{if } H\in (\frac58,\, \frac34).
  \end{array}
\right.
\end{equation*} 
\end{proposition}
\begin{proof}For simplicity, we assume that $\theta=1$. $g_1\otimes g_2$
indicates the tensor product of the two functionals $g_1$ and $g_2$, and $g^{\otimes 2}:=g\otimes g$.
It follows from the definition of double Wiener-It\^o integral and the stochastic Fubini theorem that 
\begin{align*}
   \frac{1}{2\sqrt{T}}I_2(f_T)&= \frac{1}{\sqrt{T}}\int_0^T\eta_t\,\dif B^H_t=W_T+\frac{1}{2\sqrt{T}}\big[I_2 (h_T^{\otimes 2}) \big],
\end{align*}where $h_T(u)=e^{u-T}\mathbf{1}_{[0,T]}(u)$, see also (3.4) of 
[6] or (3.21) of 
[2].  It\^o's isometry implies that 
\begin{align*}
    \E [W_T^2 ]&=\frac{1}{4 {T}}\E\big[(I_2(f_T)-I_2 (h_T^{\otimes 2}) )^2\big]=\frac{1}{2 {T}}\norm{f_T-h_T^{\otimes 2}}^2_{\FH^{\otimes 2}}\\
    &=\frac{1}{2 {T}}\left[\norm{f_T}_{\FH^{\otimes 2}}^2-2\innp{f_T,\, h_T^{\otimes 2}}_{\FH^{\otimes 2}}+\norm{h_T}_{\FH}^4\right].
\end{align*}
By Lemma~\ref{upper bound F}, there is a constant $C$ independent of $T$ such that  $\norm{h_T}\le C$. Thus, 
\begin{align*}
    \abs{\innp{f_T,\, h_T^{\otimes 2}}_{\FH^{\otimes 2}}}\le {\norm{f_T}_{\FH^{\otimes 2}}}\norm{h_T}_{\FH}^2\le C\sqrt{T}.
\end{align*}
Hence, we have 
\begin{align*}
    \abs{\E [W_T^2 ]-\sigma^2}&=\abs{\frac{1}{2 {T}}\left[\norm{f_T}_{\FH^{\otimes 2}}^2-2\innp{f_T,\, h_T^{\otimes 2}}_{\FH^{\otimes 2}}+\norm{h_T}_{\FH}^4\right]-\sigma^2 }\\
    &\le \abs{\frac{1}{2 {T}} \norm{f_T}_{\FH^{\otimes 2}}^2-\sigma^2}+\frac{C}{\sqrt{T}}\\
    &\le \frac{C}{T^{1\wedge (3-4H)}}+\frac{C}{\sqrt{T}}\le \frac{C}{T^{\gamma}},
\end{align*} where in the second inequality we use Corollary~\ref{corjjx}.
\end{proof}
The stationary solution to \eqref{OU dingyi duibi} has an integral representation
\begin{align}\label{zt pw}
     Z_t = \int_{-\infty}^t e^{-\theta(t-u)}\dif B^H_u.
\end{align}
It is known that the decay of its auto-covariance function $\rho(r)=\E[Z_rZ_0]$ is like that of a power function, see 
[7]. As a comparison, in the previous work such as 
[2; 5; 8-9], both the following expression 
\begin{align}\label{wt2 bs}
    \E[W_T^2]=\frac{2}{  {T}}\int_{[0,T]^2}\rho_1^2(t,s)\dif t\dif s,
\end{align}where $\rho_1(t,s)=\E[\eta_t\eta_s]$, and the power function decay of the function $\rho(r)$ are used to obtain the estimate \eqref{erjielieliang}, see (2.21) of 
[5] for example. 

The following identity \eqref{yigedengshi} is a corollary of \eqref{erjielieliang}, which appears in 
[2; 8] implicitly. 
\begin{corollary} \label{coro.2.3}
Let $Z_t$ be the stationary solution to \eqref{OU dingyi duibi} and $\rho(r)=\E[Z_rZ_0]$. Then 
   \begin{equation}\label{yigedengshi}
\int_0^{\infty}  \rho^2(r) \dif r=\frac14 H^2\Gamma^2(2H) (4H-1)(1-\frac{1}{\cos(2H \pi)}).
\end{equation} 
\end{corollary}
\begin{proof}Denote $V_T=\frac{1}{\sqrt {T}}\int_{0}^T \big( Z_t^2 -\E[Z_t^2] \big)\dif t$.
    It follows from (3.24) of 
    [2] that there exists a positive constant $C$ independent of $T$ such that  
\begin{align*}
    \abs{\E[W_T^2-V_T^2]}\le \frac{C}{\sqrt{T}},
\end{align*} see also the proof of Lemma 1 of 
[1]. Similar to \eqref{wt2 bs}, it is known that as $T\to \infty$
\begin{align*}
    \E[V_T^2]=\frac{2}{  {T}}\int_{[0,T]^2}\rho^2(t-s)\dif t\dif s\to 4\int_0^{\infty}  \rho^2(r) \dif r.
\end{align*} Combining the above two equation with \eqref{erjielieliang}, we obtain \eqref{yigedengshi} by means of the triangle inequality. 
    
\end{proof}

\section{Proof of Theorem~\ref{qq2b}}\label{sec.3}
In this section, we will show Theorem~\ref{qq2b} by means of two different inner product formulas according to $H\in (0,\frac12)$ and $H\in (\frac12,\frac34)$.\\
 {\bf Proof of Theorem~\ref{qq2b} when $H\in (0,\frac12)$: } 
 Denote $\beta=2H-1$. It follows from Proposition~\ref{st1-thm-original} that 
    \begin{align}
     \|f_T\|^2_{{\FH}^{\otimes2}}
    & =H^2\int_{[0,T]^4} \frac{\partial^2}{\partial t_1 \partial s_2} \left\{e^{-|t_1-s_1|-|t_2-s_2|} 1_{\{0\leq t_1, s_2\leq T\}}\right\}\notag\\
    &\times \sgn(t_2-t_1)  \abs{t_2-t_1}^{\beta}  \sgn(s_1-s_2)  \abs{s_1-s_2}^{\beta} \dif\vec{t}   \dif \vec{s}, \label{chufadian000}
    \end{align}where the partial derivative is given by:
\begin{align*}
    \frac{\partial^2}{\partial t_1 \partial s_2} \left\{e^{-|t_1-s_1|-|t_2-s_2|} 1_{\{0\leq s_2,t_1\leq T\}}\right\} &=e^{-|t_1-s_1|-|t_2-s_2|} \left[\sgn(s_1-t_1)+(\delta_0(t_1)-\delta_T(t_1)) \right]\\
    &\times \left[\sgn(t_2-s_2)+(\delta_0(s_2)-\delta_T(s_2)) \right].
\end{align*}
Substituting the above equation into \eqref{chufadian000}, we have by symmetry,
\begin{align}   \frac{1}{H^2} \|f_T\|^2_{{\FH}^{\otimes2}}= I_1(T)+2I_2(T)+I_3(T),\label{qishift2-fenjie}
\end{align} where 
\begin{eqnarray*}
    I_1(T) & = & \int_{[0,T]^4}  {e^{-|t_1-s_1|-|t_2-s_2|}} \sgn(s_1-t_1) \sgn(t_2-s_2) \\
    & & \times \ \sgn(t_2-t_1)  \abs{t_2-t_1}^{\beta}  \sgn(s_1-s_2)  \abs{s_1-s_2}^{\beta } \dif \vec{t}  \dif \vec{s} \ , \\
    I_2(T) & = & \int_{[0,T]^4}  {e^{-|t_1-s_1|-|t_2-s_2|}} \sgn(s_1-t_1) (\delta_0(s_2)-\delta_T(s_2))  \\
    & & \times \ \sgn(t_2-t_1)  \abs{t_2-t_1}^{\beta}  \sgn(s_1-s_2) \abs{s_1-s_2}^{\beta } \dif \vec{t}  \dif \vec{s} \ , \\
    I_3(T) & = & \int_{[0,T]^4}  {e^{-|t_1-s_1|-|t_2-s_2|}} (\delta_0(t_1)-\delta_T(t_1)) (\delta_0(s_2)-\delta_T(s_2))  \\
    & & \times \  \sgn(t_2-t_1)  \abs{t_2-t_1}^{\beta}  \sgn(s_1-s_2)  \abs{s_1-s_2}^{\beta } \dif \vec{t}  \dif \vec{s} \ .
\end{eqnarray*}
We will compute these three items one by one. First, by symmetry and Lemma~\ref{L1T calculat}, we have that
\begin{align}
    I_1(T)& =4\int_{0}^T e^{-t_1}\dif t_1 \int_{[0,t_1]^3}  {e^{ s_1-|t_2-s_2|}} (t_1-t_2)^{\beta}   \abs{s_1-s_2}^{\beta} \sgn(t_2-s_2) \sgn(s_1-s_2)\dif t_2 \dif \vec{s}\notag\\
    &=2 (4H-1) \Gamma^2(2H)\left[1-\frac{ 1}{\cos(2H\pi)}\right]\Big[T-\frac{4H+1}{2} \Big]\notag\\
    &-8H^2\Gamma^2(2H)+\frac{2}{4H-3}T^{4H-2} +O(T^{4H-3}).\label{i1t jielun-0000}
\end{align}

Second, we compute the term $I_2(T)$. 
By the definition of Dirac delta function, we
have 
\begin{align}
    I_2(T) & = \int_{[0,T]^3}
     {e^{-|t_1-s_1|-t_2}} \sgn(s_1-t_1)   \sgn(t_2-t_1)  \abs{t_2-t_1}^{\beta} s_1^{\beta } \dif \vec{t}  \dif s_1\notag \\
     &+ \int_{[0,T]^3}
     {e^{-|t_1-s_1|-(T-t_2)}} \sgn(s_1-t_1)   \sgn(t_2-t_1)  \abs{t_2-t_1}^{\beta} (T-s_1)^{\beta } \dif \vec{t}  \dif s_1.
\end{align}Making the change of variable $s_1'=T-s_1,\, t_1'=T-t_1,\, t_2'=T-t_2$, we have 
\begin{align*}
   & \int_{[0,T]^3}
     {e^{-|t_1-s_1|-(T-t_2)}} \sgn(s_1-t_1)   \sgn(t_2-t_1)  \abs{t_2-t_1}^{\beta} (T-s_1)^{\beta } \dif \vec{t}  \dif s_1\\
   &=  \int_{[0,T]^3}
     {e^{-|t_1'-s_1'|-t_2'}} \sgn(s_1'-t_1')   \sgn(t_2'-t_1')  \abs{t_2'-t_1'}^{\beta} (s_1')^{\beta } \dif \vec{t'}  \dif s_1'
\end{align*}
Hence, we can simplify $$I_2(T) = 2\int_{[0,T]^3}
     {e^{-|t_1-s_1|-t_2}} \sgn(s_1-t_1)   \sgn(t_2-t_1)  \abs{t_2-t_1}^{\beta} s_1^{\beta } \dif \vec{t}  \dif s_1. $$
We can split $[0, T]^3$ into the sub-regions:
\begin{align*}
   \bar{\Delta}_1&=\set{0\le t_1\le s_1 \le t_2\le T}\cup \set{0\le t_2\le s_1 \le t_1\le T},  \\
   \bar{\Delta}_2&=\set{0\le t_1\le t_2\le s_1 \le T}\cup \set{0\le t_2 \le t_1\le s_1\le T},\\
   \bar{\Delta}_3&=\set{0\le s_1 \le t_1\le t_2\le  T}\cup \set{0\le s_1 \le t_2 \le t_1\le T}.
\end{align*}
The integrations of the integrand of $I_2(T)$ over the sub-regions are denoted as $I_{21}(T),\,I_{22}(T)$ and $I_{23}(T)$ respectively. By symmetry, $I_{23}(T)=0$. Then we have 
\begin{align}
 I_2(T)=2\left[I_{21}(T)+I_{22}(T)+I_{23}(T)\right]=2\left[I_{21}(T)+I_{22}(T) \right].\label{I2T bds} 
\end{align}
Making the change of variable $x=\abs{t_1-t_2}, y=s_1-(t_1\wedge t_2)$, we have 
\begin{align}
 I_{21}(T)&= \int_{0\le y\le x,s_1\le T,\,s_1+x-y\le T}e^{-(x+s_1)}x^{\beta} s_1^{\beta} \left(1+e^{2y}\right) \dif x\dif s_1\dif y\notag\\
 &= \int_{[0,T]^2}e^{-(x+s)}x^{\beta}s^{\beta} \dif x\dif s \notag\\
 & \quad \times \left( x\wedge s- 0\vee(s+x-T)+\frac12 e^{2( x\wedge s)}- \frac12 e^{2 \left(0\vee(s+x-T)\right)}\right).\label{i21 fenjiezhk}
\end{align}
It is evidence for any $\gamma>0$,
\begin{align}
   &- \int_{[0,T]^2}e^{-(x+s)}x^{\beta}s^{\beta} \left( 0\vee(s+x-T) \right)\dif x\dif s \notag\\
   &=\int_{[0,T]^2,s+x\ge T}e^{-(x+s)}x^{\beta}s^{\beta} (T-x-s) \dif x\dif s=O(\frac{1}{T^{\gamma}}).\label{j21-1}
\end{align}
By Lemma~\ref{asymptotic expansion key-001}, we have that 
\begin{align}
&-\frac12 \int_{[0,T]^2}e^{-(x+s)}x^{\beta}s^{\beta}  e^{2 \left(0\vee(s+x-T)\right)}\dif x\dif s  \notag\\
&=-\frac12 \int_{[0,T]^2}e^{-(x+s)}x^{\beta}s^{\beta}   \dif x\dif s+\frac12\int_{[0,T]^2,x+s\ge T}e^{-(x+s)}x^{\beta}s^{\beta}   \dif x\dif s\notag\\
&-\frac12 \int_{[0,T]^2,x+s\ge T}e^{x+s-2T}x^{\beta}s^{\beta}   \dif x\dif s\notag\\
&=-\frac12\Gamma^2(2H)-\frac12\left[T^{2\beta}-2\beta T^{2\beta-1} +O(T^{2\beta-2})\right].\label{j21-2}
\end{align}
The symmetry and the integration by parts imply that for any $\gamma>0$,
\begin{align}
   & \int_{[0,T]^2}e^{-(x+s)}x^{\beta}s^{\beta} \left( x\wedge s +\frac12 e^{2( x\wedge s)} \right)\dif x\dif s \notag\\
   &=2\int_0^{T} e^{-s}s^{\beta} \dif s\int_0^s x^{1+\beta}\dif (-e^{-x})+ \int_{0\le x\le s\le T} e^{x-s}x^{\beta}s^{\beta}\dif x\dif s \notag\\
   &=-2\int_0^{T} e^{-2s}s^{2\beta+1}\dif s+(\beta+1)\int_{[0,T]^2}e^{-x-s}x^{\beta}s^{\beta}\dif x\dif s+ \int_{0\le x\le s\le T} e^{x-s}x^{\beta}s^{\beta}\dif x\dif s \notag\\
   &=-2^{1-4H}\Gamma(4H)+2H\Gamma^2(2H)+ \int_{0\le x\le s\le T} e^{x-s}x^{\beta}s^{\beta}\dif x\dif s+O(\frac{1}{T^{\gamma}}).\label{j21-3}
\end{align}  Substituting \eqref{j21-1}-\eqref{j21-3} into \eqref{i21 fenjiezhk}, we have 
\begin{align}
    I_{21}(T)&=-2^{1-4H}\Gamma(4H)+(2H-\frac12)\Gamma^2(2H)+ \int_{0\le x\le s\le T} e^{x-s}x^{\beta}s^{\beta}\dif x\dif s\notag\\
    &-\frac12\left[T^{2\beta}-2\beta T^{2\beta-1} +O(T^{2\beta-2})\right].\label{j21-bds}
\end{align}
Making again the change of variable $x=\abs{t_1-t_2}, y=s_1-(t_1\wedge t_2)$  and using the integration by parts, we have that for any $\gamma>0$
\begin{align}
 I_{22}(T)&= \int_{0\le x\le y\le s_1\le T}  \left(e^{-x-s_1}-e^{x-s_1}\right) x^{\beta}s_1^{\beta}\dif x\dif y\dif s_1\notag\\
 &=\int_{0\le x\le s\le T} e^{x-s}s^{\beta}x^{\beta} (x-s)\dif x\dif s +2\int_0^Te^{-2s}s^{2\beta+1}\dif s-e^{-T}T^{1+\beta}\int_0^Te^{-x}x^{\beta}\dif x\notag\\
 &=\int_{0\le x\le s\le T} e^{x-s}s^{\beta}x^{\beta} (x-s)\dif x\dif s +2^{1-4H}\Gamma(4H)+O(\frac{1}{T^{\gamma}}).\label{i22 fenjiezhk}
\end{align}The integration by parts and Lemma~\ref{asymptotic expansion key-0} imply for any $\gamma>0$
\begin{align}
    &\int_{0\le x\le s\le T} e^{x-s}x^{\beta}s^{\beta} \dif x\dif s +\int_{0\le x\le s\le T} e^{x-s}s^{\beta}x^{\beta} (x-s)\dif x\dif s\notag\\
    &=-(2\beta+1)\int_{0\le x\le s\le T} e^{x-s}x^{\beta}s^{\beta} \dif x\dif s+T^{1+\beta}e^{-T}\int_0^T e^x x^{\beta}\dif x\notag\\
    &=-(2\beta+1)\int_{0\le x\le s\le T} e^{x-s}x^{\beta}s^{\beta} \dif x\dif s+T^{1+\beta}\left[T^{\beta}-\beta T^{\beta-1} +\beta (\beta-1)T^{\beta-2}+O(T^{\beta-3})\right].\label{zhongjianjieg}
\end{align}

Substituting \eqref{j21-bds} and \eqref{i22 fenjiezhk} and \eqref{zhongjianjieg} into \eqref{I2T bds}, we obtain  
\begin{align}
    I_2(T)&=2\left[\int_{0\le x\le s\le T} e^{x-s}x^{\beta}s^{\beta} \dif x\dif s +\int_{0\le x\le s\le T} e^{x-s}s^{\beta}x^{\beta} (x-s)\dif x\dif s\right]\notag\\
    &+(4H-1)\Gamma^2(2H)-\left[T^{2\beta}-2\beta T^{2\beta-1} +O(T^{2\beta-2})\right]\notag\\
    &=-2(2\beta+1)\int_{0\le x\le s\le T} e^{x-s}x^{\beta}s^{\beta} \dif x\dif s +(4H-1)\Gamma^2(2H)\notag\\
   &+2T^{2\beta+1}-(2\beta+1)T^{2\beta}+2\beta^2 T^{2\beta-1}+O(T^{2\beta-2})\notag\\
    &=(4H-1)\Gamma^2(2H)\big[1-\frac{1}{\cos (2H\pi)}\big]+\frac{2\beta}{2\beta-1}T^{2\beta-1}+O(T^{2\beta-2}),\label{i2t jielun}
\end{align} where in the last line we use Lemma~\ref{asymptotic expansion key}.
    
Using the definition of Dirac delta function and the basic identities 
\begin{align*}
   \int_0^Te^{-u}(T-u)^{\beta}\dif u&=e^{-T } \int_0^T e^x x^{\beta} \dif x=T^{\beta} -\beta T^{\beta-1} +O(T^{\beta-2}), \\
   \int_0^Te^{u-T}(T-u)^{\beta}\dif u&=\int_0^T e^{-x} x^{\beta} \dif x =\Gamma(2H)+O(\frac{1}{T^{\gamma}}),
\end{align*} for any $\gamma>0$,
it is evident that 
\begin{align}
    I_3(T)
    &=\int_{[0,T]^4}  {e^{-|t_1-s_1|-|t_2-s_2|}} (\delta_0(t_1)-\delta_T(t_1)) (\delta_0(s_2)-\delta_T(s_2))  \notag\\
    & \quad \times \sgn(t_2-t_1)  \abs{t_2-t_1}^{\beta}  \sgn(s_1-s_2)  \abs{s_1-s_2}^{\beta } \dif \vec{t}  \dif \vec{s}\notag\\
    &=2\Gamma^2(2H)+{2}T^{2\beta}-4\beta T^{2\beta-1}+O(T^{2\beta-2}).\label{i3t jielun}
    \end{align}
Finally, inserting the computation results \eqref{i1t jielun-0000}, \eqref{i2t jielun} and \eqref{i3t jielun} into \eqref{qishift2-fenjie}, we obtain the desired results.

 {\bf Proof of Theorem~\ref{qq2b} when $H\in [\frac12,1)$: } We will
discuss exclusively the case $H\in (\frac12,1)$ since the case $H =\frac12$ is easy. Denote $\beta=2H-2$. It\^o's isometry implies that 
\begin{align}
     \|f_T\|^2_{{\FH}^{\otimes2}}
    & =H^2(2H-1)^2 \int_{[0,T]^4}  e^{-\abs{t_1-s_1}-\abs{t_2-s_2}}  \abs{t_2-t_1}^{\beta}  \abs{s_1-s_2}^{\beta} \dif\vec{t}   \dif \vec{s}. \label{chufadian100}
\end{align}The symmetry implies that 
\begin{align}
    \frac{1}{H^2} \|f_T\|^2_{{\FH}^{\otimes2}}
    & =4(2H-1)^2 \int_{0}^T e^{-t_1}\dif t_1 \int_{[0,t_1]^3}  {e^{ s_1-|t_2-s_2|}} (t_1-t_2)^{\beta}   \abs{s_1-s_2}^{\beta}  \dif t_2 \dif \vec{s}. \label{chufadian100-1}
\end{align}Using Lemma~\ref{L2T calculat}, we have that 
\begin{align*}
 \frac{1}{H^2} \|f_T\|^2_{{\FH}^{\otimes2}}
    & =4(1+\beta)^2\Bigg\{\frac12\Gamma^2({1+\beta})\left[(2\beta+3)(1-\frac{1}{\cos(\beta \pi)})(T-\frac{2\beta+3}{2})-(\beta+1)(\beta+3) \right]\notag\\
    &+\frac{1}{(\beta+1)(2\beta+1)}T^{2\beta+2}- \frac{2}{2\beta+1}T^{2\beta+1}+O(T^{2\beta})\Bigg\}\notag\\
    &=2\Gamma^2(2H)\left[(4H-1)(1-\frac{1}{\cos(2H \pi)})(T-\frac{4H-1}{2})-(2H-1)(2H+1) \right]\notag\\
    &+\frac{4(2H-1)}{4H-3}T^{4H-2}-\frac{8(2H-1)^2}{4H-3}T^{4H-3}+O(T^{4H-4}). 
\end{align*}{\hfill\large{$\Box$}}  
\section{Appendix}
For a positive function $\phi$, we say a real-valued function $f$ satisfies $f=o(\phi)$ if $\frac{f(u)}{\phi(u)}\to 0$ as $u\to \infty$. Especially, the notation $h(u)=o(1)$ means that $h(u)\to 0$ as $u\to \infty$.  Lemma~\ref{upper bound F} is well known, see 
[1] for example.

\begin{lemma} \label{upper bound F}
 Assume $\beta>-1$.  Then there exists a constant $C>0$ such that for any  $s\in [0,\infty)$,
\begin{align*}
e^{-  s}\int_0^{s} e^{  r} r^{\beta }\dif r&\le C \times\big(s^{\beta+1}\mathbbm{1}_{[0,1]}(s) + s^{\beta}\mathbbm{1}_{ (1,\,\infty)}(s)\big).
\end{align*}
Especially, when $\beta\in (-1,0)$, there exists a constant $C>0$ such that for any  $s\in [0,\infty)$,
\begin{align*}
e^{-  s}\int_0^{s} e^{  r} r^{\beta }\dif r\le C \times(1\wedge s^{\beta}).
\end{align*}
\end{lemma}

\begin{lemma}\label{asymptotic expansion key-0} 
    Suppose  $\beta\in(-1,0)$. The following asymptotic expansion
    \begin{align*}  
    e^{-T } \int_0^T e^x x^{\beta} \dif x =T^{\beta}-\beta T^{\beta-1} +\beta (\beta-1)T^{\beta-2}+O(T^{\beta-3})
    \end{align*}
     holds as $T\to \infty$.
\end{lemma}
\begin{proof}
    It can be obtain by a slight change of proof of Lemma A.3 of 
    [1]. 
\end{proof}
\begin{lemma}\label{asymptotic expansion key-001} 
    Suppose  $\beta\in(-1,0)$. The following asymptotic expansion
    \begin{align*}  
    \int_{[0, T]^2,x+z\ge T }  e^{z +x-2T}   x^{\beta}  z ^{\beta}\dif x \dif z=T^{2\beta}-2\beta T^{2\beta-1} +  {\beta}(3\beta-2)T^{2\beta-2}+O(T^{2\beta-3}).
    \end{align*}
     holds as $T\to \infty$.
\end{lemma}
\begin{proof}
 Let the function $g(x)=e^{2T}T^{2\beta-2}$. Since  $\lim_{x\to \infty}g(x)= \infty$ and
$g'(x)\neq 0 $ in the neighborhood of $\infty$, then the lemma can  be proved by applying L'H\^{o}pital's rule to show the limit
    \begin{align*}
       & \lim_{T\to \infty} \frac{   \int_{[0, T]^2,x+z\ge T }  e^{z +x-2T}   x^{\beta}  z ^{\beta}\dif x \dif z-(T^{2\beta}-2\beta T^{2\beta-1} +{\beta}(3\beta-2)T^{2\beta-2})}{T^{2\beta-3}}\\
        &=\lim_{T\to \infty} \frac{\int_0^T e^x x^{\beta} \dif x\int_{T-x}^T e^z z^{\beta}\dif z- e^{2T}(T^{2\beta}-2\beta T^{2\beta-1}+{\beta}(3\beta-2)T^{2\beta-2})}{e^{2T}T^{2\beta-3}} \\
        &= {\beta}(\beta-1)(3\beta-2) .
    \end{align*}
\end{proof}
\begin{lemma}\label{linjiezhi}
The following asymptotic expansion holds as $T\to \infty$:
\begin{align}\label{IT zk}
  I(T):=    \int_0^1\frac{1}{\sqrt{1-t}}\frac{1-e^{-tT}}{t} \dif t= \log T +2\log 2+\gamma-\frac{1}{2T}-\frac{3}{8T^2}+O(T^{-3}),
\end{align}where $\gamma=\int_0^1\frac{1-e^{-u}-e^{-1/u}}{u}\dif u$ is the Euler’s constant.
\end{lemma}
\begin{proof}
First, we rewrite 
\begin{align}\label{qishi IT}
    I(T)=\left(\int_0^{\frac{1}{T}}+\int_{\frac{1}{T}}^1\right)\frac{1}{\sqrt{1-t}}\frac{1-e^{-tT}}{t} \dif t:=I_1(T)+I_2(T).
\end{align}   
The change of variable $u=tT$ and the asymptotic expansion \begin{align}\label{erxiangshizk}
    (1+x)^{\alpha}=1+\alpha x+\frac{\alpha(\alpha-1)}{2}x^2+O(x^3),\qquad \text{as \,\,} x\to 0
\end{align} imply that 
\begin{align}
    I_1(T)&=\int_0^1\frac{1-e^{-u}}{u}\frac{1}{\sqrt{1-\frac{u}{T}}}\dif u\notag\\
    &=\int_0^1\frac{1-e^{-u}}{u} \left(1+\frac{1}{2}\frac{u}{T}+\frac{3}{8}(\frac{u}{T})^2+O((\frac{u}{T})^3)\right)\dif u\notag\\
    &=\int_0^1\frac{1-e^{-u}}{u}  \dif u +\frac{1}{2eT}+\frac{3}{8T^2}(2e^{-1}-\frac12)+O(T^{-3}).\label{I1T zhk}
\end{align}The change of variable $u=\sqrt{1-t}$ implies that $$\int_0^1[\frac{1}{t\sqrt{1-t}}-\frac{1}{t}]\dif t=\int_0^1\frac{1}{\sqrt{1-t}(1+\sqrt{1-t})} \dif t =2\int_0^1\frac{1}{1+u}\dif u=2\log 2.$$
Hence, we have 
\begin{align}
    I_2(T)&=\log T+\int_{\frac{1}{T}}^1[\frac{1}{t\sqrt{1-t}}-\frac{1}{t}]\dif t-\int_{\frac{1}{T}}^1\frac{e^{-tT}}{t\sqrt{1-t}} \dif t\notag\\
    &=\log T+2\log 2 - \int_0^{\frac{1}{T}}\frac{1}{t}[\frac{1}{\sqrt{1-t}}-1]\dif t -\int_{\frac{1}{T}}^1\frac{e^{-tT}}{t\sqrt{1-t}} \dif t.\label{I2T zhk}
\end{align}
The asymptotic expansion \eqref{erxiangshizk} implies that 
\begin{align}
    \int_0^{\frac{1}{T}}\frac{1}{t}[\frac{1}{\sqrt{1-t}}-1]\dif t&=\int_0^{\frac{1}{T}}\frac{1}{t}[\frac12 t+\frac{3}{8}t^2+\frac{5}{16}t^3+O(t^4)]\dif t\notag\\
    &=\frac{1}{2T}+\frac{3}{16T^2}+O(T^{-3}).\label{I2T zhj 1 zk}
\end{align}The change of variable $u=\frac{1}{tT}$ and the asymptotic expansion \eqref{erxiangshizk} imply that 
\begin{align}
    \int_{\frac{1}{T}}^1\frac{e^{-tT}}{t\sqrt{1-t}} \dif t&=\int_{\frac{1}{T}}^1 \frac{e^{-\frac{1}{u}}}{u}\left(1-\frac{1}{Tu}\right)^{-\frac12}\dif u \notag\\
    &=\int_{\frac{1}{T}}^1 \frac{e^{-\frac{1}{u}}}{u}\left(1+ \frac{1}{2Tu}+\frac{3}{8}(\frac{1}{Tu})^2 +O((\frac{1}{Tu})^3) \right)\dif u\notag\\
    &=\int_0^1\frac{e^{-\frac{1}{u}}}{u}\dif u +\frac{1}{2eT}+\frac{3}{4eT^2}+O(T^{-3}). \label{I2T zhj 2 zk}
\end{align}
Substituting the asymptotic expansions \eqref{I1T zhk}-\eqref{I2T zhj 2 zk} into \eqref{qishi IT}, we obtain \eqref{IT zk}.
\end{proof}

\begin{lemma}\label{asymptotic expansion key}
    Suppose $\beta\in (-1,\frac12)$ 
    The following asymptotic expansion
    \begin{align*}
  &  \int_{0\le x \le z \le T } e^{x-z} x^{\beta}  z ^{\beta}\dif x \dif z \\
  &  =\left\{
      \begin{array}{ll}
\log T +2\log 2+\gamma-\frac{1}{2T}-\frac{3}{8T^2}+O(T^{-3}), & \,\, \beta= -\frac12,\\
-\frac{\Gamma^2(1+\beta)}{2\cos(\beta\pi)} +\frac{T^{\delta}}{\delta} {-}\frac12T^{\delta-1}+ \frac{\beta(\beta-1)}{\delta-2}T^{\delta-2}{-}\frac{\beta(\beta-2)}{2}T^{\delta-3}+O(T^{\delta-4}),& \,\, \beta\neq -\frac12, 
 \end{array}
\right.     
    \end{align*}
     holds, where ${\delta}:=2\beta+1$ and $\gamma$ 
     is the Euler’s constant.
\end{lemma}  
\begin{proof}The case of $\beta=-\frac12$ is from the change of variable $x=(1-t)z$ and Lemma~\ref{linjiezhi}.
The case $\beta\neq -\frac12$ is a special case of Lemma~A.3 of 
    [1]. In fact,  it follows from Lemma~A.3 of 
    [1] that  when $\beta\in (-1,0)$, the following
asymptotic expansion holds.
       \begin{align*}
  &  \int_{0\le x \le z \le T } e^{x-z} x^{\beta}  z ^{\beta}\dif x \dif z \\
  &  =\left\{
      \begin{array}{ll}
\Gamma(\delta+1)\mathrm{B}(1+\beta, -\delta)+\delta^{-1}T^{\delta} -\frac12T^{\delta-1}+O(T^{\delta-2}), & \quad  \delta\in (-1, 0),\\
\log T +o(\log T), & \quad \delta= 0,\\
\delta^{-1}T^{\delta}+\beta\Gamma(\delta)\mathrm{B}(1+\beta,1-\delta){-}\frac12T^{\delta-1}+O(T^{\delta-2}),& \quad \delta\in (0,1),
 \end{array}
\right.     
    \end{align*}
     
Using the basic relations $\Gamma(1+z)=z\Gamma(z) $ and $\Gamma(z)\Gamma(1-z)=\frac{\pi}{\sin(\pi z)}$ and $B(a,b)=\frac{\Gamma(a)\Gamma(b)}{\Gamma(a+b)}$, we have that when $\beta\in (-1,-\frac12)$, i.e., $\delta\in (-1,0)$,
\begin{align*}
    \Gamma(\delta+1)\mathrm{B}(1+\beta, -\delta)&=\Gamma^2(1+\beta)\frac{\Gamma(2\beta+2)\Gamma(-1-2\beta)}{\Gamma(1+\beta)\Gamma(-\beta)}\\
    &=\Gamma^2(1+\beta)\frac{\sin(\beta\pi)}{\sin[(1+2\beta)\pi]}=\frac{-1}{2\cos(\beta\pi)}\Gamma^2(1+\beta),
\end{align*} and when $\beta\in (-\frac12,0)$, i.e., $\delta\in (0,1)$,
\begin{align*}
   \beta\Gamma(\delta)\mathrm{B}(1+\beta,1-\delta)&=\Gamma^2(1+\beta)\frac{\beta\Gamma(\delta)\Gamma(1-\delta)}{\Gamma(1+\beta)\Gamma(1-\beta)}\\
    &=\Gamma^2(1+\beta)\frac{\sin(\beta\pi)}{\sin\big[(2\beta+1)\pi\big]}=\frac{-1}{2\cos(\beta\pi)}\Gamma^2(1+\beta).
\end{align*}
When $\beta\in [0,\frac12)$, i.e., $\delta\in[1,2)$, we can repeat the whole proof of Lemma~A.3 of 
[1] for the case of $\delta\in(0,1)$.

Finally, the addition terms $\frac{\beta(\beta-1)}{\delta-2}T^{\delta-2}-{\frac{\beta(\beta-2)}{2}}T^{\delta-3}
$ can be obtain by a  slightly routine change of proof of Lemma~A.3 of 
[1].
\end{proof}

\begin{lemma}\label{L1T calculat}
Suppose $H\in (0,\frac34)$ and $\beta=2H-1$.
\begin{align}
  L_1(T)&:=\int_{0}^T e^{-t_1}\dif t_1 \int_{[0,t_1]^3}  {e^{ s_1-|t_2-s_2|}} (t_1-t_2)^{\beta}   \abs{s_1-s_2}^{\beta} \sgn(t_2-s_2) \sgn(s_1-s_2)\dif t_2 \dif \vec{s}\notag\\    
  &=\frac12  \Gamma^2(2H)\left[(4H-1)-\frac{4H-1}{\cos(2H\pi)}\right]\Big[T-\frac{4H+1}{2} \Big]\notag\\
    &-2H^2\Gamma^2(2H)+\frac{1}{2(4H-3)}T^{4H-2} -\frac{4H-2}{4H-3}T^{4H-3}+O(T^{4H-4}), \label{fenjie 001}
\end{align}where the term $\frac{4H-1}{-\cos(2H\pi)}$ in \eqref{fenjie 001} is illustrated as its limit $\frac{2}{\pi}$ when $H=  \frac14$.
\end{lemma}
\begin{proof}

The domain $[0,T]\times [0,t_1]^3$ can be divided into six disjoint regions according
to the distinct orders of $t_2,s_1,s_2$: 
\begin{equation}\label{domain divided}
\begin{aligned}
\Delta_1=\set{0\le t_2\le s_2\le s_1\le t_1 \le T }, \quad \Delta_{1'}=\set{0\le t_2\le  s_1\le s_2\le t_1 \le T },\\
\Delta_2=\set{0\le s_2\le t_2 \le s_1\le t_1 \le T },  \quad\Delta_{2'}=\set{0\le  s_1\le t_2\le s_2\le t_1 \le T },\\
\Delta_3=\set{0\le s_2\le s_1\le t_2\le t_1 \le T },  \quad\Delta_{3'}=\set{0\le  s_1\le s_2\le t_2\le t_1 \le T }.
\end{aligned}
\end{equation}
For $i=1,2,3,1',2',3'$, denote
$$L_{1i}(T)=\int_{\Delta_i}  {e^{ s_1-t_1-|t_2-s_2|}} (t_1-t_2)^{\beta}   \abs{s_1-s_2}^{\beta} \sgn(t_2-s_2) \sgn(s_1-s_2)\dif \vec{t} \dif \vec{s},$$
so $L_1(T) = \sum_{i}L_{1i}(T)$. Making substitution $x=\abs{s_1-s_2}, y=(s_1\vee s_2)-t_2, z=t_1-t_2$, we have that
\begin{align}
   L_{11}(T)+ L_{11'}(T)&= \int_{0\le x\le y\le z\le t_1\le T }  \left({e^{-z -x}-e^{x-z}}\right) x^{\beta}  z ^{\beta} \dif x \dif y \dif z\dif  t_1\notag\\
   &=\int_{0\le x \le z \le T }  \left({e^{-z -x}-e^{x-z}}\right) x^{\beta}  z ^{\beta} (z-x)(T-z) \dif x \dif z,\label{diyibufen-01}
\end{align}
and that 
\begin{align}
   &L_{12}(T)+ L_{12'}(T)\notag\\
   &= \int_{0\le y\le x,z\le x+z-y\le t_1\le T }  \left({e^{-z -x}+e^{-z-x+2y}}\right) x^{\beta}  z ^{\beta} \dif x \dif y \dif z\dif  t_1\notag\\
   &=\int_{[0, T]^2 }  e^{-z -x}  x^{\beta}  z ^{\beta}    \dif x \dif z\int^{x\wedge z}_{0\vee \left(x+z-T \right)}\left(1+e^{2y}\right)\left(T -\left( x+z-y\right)\right) \dif y\notag\\
   &=\frac14 \int_{[0, T]^2,x+z\ge T }  e^{z +x-2T}  x^{\beta}  z ^{\beta}\dif x \dif z + \frac12\int_{[0, T]^2 }  e^{-z -x}  x^{\beta}  z ^{\beta}  (T-x\vee z)^2  \dif x \dif z \notag\\
   & \quad -\frac12 \int_{[0, T]^2,x+z\le T }  e^{-z -x}  x^{\beta}  z ^{\beta} \left[ (T-x- z)^2 +(T-x-z)-\frac12 \right ]\dif x \dif z\notag\\
   & \quad +\frac12 \int_{[0, T]^2 }  e^{(x\wedge z)-(x\vee z)}  x^{\beta}  z ^{\beta}  (T-x\vee z -\frac12) \dif x \dif z , \label{dierbufen}
\end{align}
{ where in the last line, we firstly calculate the integral $\int^{x\wedge z}_{0\vee \left(x+z-T \right)}\left(T -\left( x+z-y\right)\right) \dif y $ directly and then use integration by parts to 
$\frac12\int^{x\wedge z}_{0\vee \left(x+z-T \right)} \left(T -\left( x+z-y\right)\right) \dif e^{2y}$, and finally rearrange all the terms according to $x+z\ge T$ and $x+z\le T$.}
\\
It is straightforward to see that
\begin{align}
   L_{13}(T)+ L_{13'}(T)= \int_{\Delta_3\cup\Delta_{3'} }  {e^{ s_1-t_1-t_2+s_2}} (t_1-t_2)^{\beta}   \abs{s_1-s_2}^{\beta} \sgn(s_1-s_2)\dif \vec{t} \dif \vec{s}=0. \label{disanbufen}
\end{align}
Substituting the above identities \eqref{diyibufen-01}-\eqref{disanbufen} into  $L_1(T) = \sum_{i}L_{1i}(T)$, we have 
\begin{align}\label{i1Tfenjie}
    L_1(T)=J_1(T)+J_2(T),
\end{align}
where 
\begin{align}
    J_2(T) &=  \frac12 \int_{[0, T]^2 }  e^{(x\wedge z)-(x\vee z)}  x^{\beta}  z ^{\beta}  (T-x\vee z -\frac12) \dif x \dif z \notag\\
    &  -\int_{0\le x \le z \le T } e^{x-z} x^{\beta}  z ^{\beta} (z-x)(T-z) \dif x \dif z, \label{j2 biaodashi}\\
    J_1(T) & =  \frac14 \int_{[0, T]^2,x+z\ge T }  e^{z +x-2T}  x^{\beta}  z ^{\beta}\dif x \dif z+\frac12\int_{[0, T]^2 }  e^{-z -x}  x^{\beta}  z ^{\beta}  (T-x\vee z)^2  \dif x \dif z \notag\\
    &  - \ \frac12 \int_{[0, T]^2,x+z\le T }  e^{-z -x}  x^{\beta}  z ^{\beta} \left[ (T-x- z)^2 +(T-x-z)-\frac12 \right ]\dif x \dif z\notag\\
    &  +\int_{0\le x \le z \le T }   e^{-z -x}  x^{\beta}  z ^{\beta} (z-x)(T-z) \dif x \dif z\notag\\
   &=\frac12 \int_{[0, T]^2 }   e^{-z -x}  x^{\beta}  z ^{\beta}  \left[ (T-x\vee z)^2  - (T-x- z)^2 -(T-x-z)+\frac12 \right ]\dif x \dif z \notag\\
 & +\int_{0\le x \le z \le T }   e^{-z -x}  x^{\beta}  z ^{\beta} (z-x)(T-z) \dif x \dif z +\frac14[T^{2\beta}-2\beta T^{2\beta-1}  +  {\beta}(3\beta-2)T^{2\beta-2}+O(T^{2\beta-3})],\label{j1 biaodashi} 
\end{align}
  where in the last line, we use Lemma~\ref{asymptotic expansion key-001} to obtain that when $H\in (0,\frac12)$,  
\begin{align*}
& \int_{[0, T]^2,x+z\ge T } \left( e^{z +x-2T} +2 e^{-z -x}   \left[ (T-x- z)^2 +(T-x-z)-\frac12 \right ] \right) x^{\beta}  z ^{\beta}\dif x \dif z\notag\\
&=\int_{[0, T]^2,x+z\ge T } e^{z +x-2T}  x^{\beta}  z ^{\beta}\dif x \dif z+O(T^{2\beta-2})=T^{2\beta}-2\beta T^{2\beta-1}  +  {\beta}(3\beta-2)T^{2\beta-2}+O(T^{2\beta-3}).
\end{align*} 

 According to the difference of squares property $a^b-b^2=(a+b)(a-b)$, 
we yield  
\begin{align}
   & \frac12\int_{[0, T]^2 }   e^{-z -x}  x^{\beta}  z ^{\beta}  \left[ (T-x\vee z)^2  - (T-x- z)^2 -(T-x-z)+\frac12 \right ]\dif x \dif z\notag\\
   &=T\int_{[0, T]^2 } e^{-z -x}  x^{\beta}  z ^{\beta}[ (x\wedge z) -\frac12]\dif x \dif z \notag\\
   &-\frac12\int_{[0, T]^2  } e^{-z -x}  x^{\beta}  z ^{\beta} \left[(x\wedge z) \left(x\wedge z+2(x\vee z)\right) -x-z-\frac12\right].\label{bigg01}
\end{align}
 The symmetry implies that  for any $\gamma>0$, 
\begin{align}
  & \frac12 \int_{[0, T]^2  } e^{-z -x}  x^{\beta}  z ^{\beta} \left[(x\wedge z)  \left(x\wedge z+2(x\vee z)\right) -x-z-\frac12\right]\dif x\dif z\notag\\
  &=\int_{0\le x\le z\le T}e^{-z -x}  x^{2+\beta}  z ^{\beta} \dif x\dif z+\int_{[0, T]^2  } e^{-z -x}  x^{1+\beta}  z ^{\beta}(z-1)\dif x \dif z\notag\\
  &-\frac14\int_{[0, T]^2  } e^{-z -x}  x^{\beta}  z ^{\beta}\dif x \dif z\notag\\
  &=\int_{0\le x\le z\le T}e^{-z -x}  x^{2+\beta}  z ^{\beta} \dif x\dif z+(\beta(1+\beta)-\frac14)\Gamma^2(1+\beta) +o(\frac{1}{T^{\gamma}}),\label{di0ge}
  \end{align}The integration by parts implies for any $\gamma>0$, 
  \begin{align}
     & \int_{0\le x\le z\le T}e^{-z -x}  x^{2+\beta}  z ^{\beta} \dif x\dif z=\int_{0}^T e^{-z}z^{\beta}\dif z\int_0^{z} x^{2+\beta}\dif (-e^{-x})\notag\\
      &=-\int_0^T e^{-2z}[z^{2\beta+2}+(2+\beta)z^{2\beta+1}]\dif z+(\beta+2)(\beta+1)\int_{0\le x\le z\le T}e^{-x-z} x^{\beta}  z ^{\beta} \dif x\dif z\notag\\
      &=-(2\beta+3)\Gamma(2\beta+2)2^{-(2+2\beta))}+\frac12(1+\beta)(\beta+2)\Gamma^2(1+\beta) +o(\frac{1}{T^{\gamma}}),
      \label{dierge}
  \end{align}
Substituting \eqref{dierge} into \eqref{di0ge}, we have for any $\gamma>0$,
 \begin{align} 
  & \frac12 \int_{[0, T]^2  } e^{-z -x}  x^{\beta}  z ^{\beta} \left[(x\wedge z)  \left(x\wedge z+2(x\vee z)\right) -x-z-\frac12\right]\dif x\dif z\notag\\&=\left((\beta+1)(\frac32\beta+1)-\frac14\right)\Gamma^2(1+\beta)-(2\beta+3)\Gamma(2\beta+2)2^{-(2+2\beta))}+o(\frac{1}{T^{\gamma}}).\label{diyibufen}
\end{align} In the same way, we have  
\begin{align}
  & \int_{0\le x \le z \le T }   e^{-z -x}  x^{\beta}  z ^{1+\beta} (z-x) \dif x \dif z\notag\\  
  &=\int_0^Te^{-x}x^{\beta}\dif x\int_x^T z^{\beta+2}\dif (-e^{-z}) -\frac12\Gamma^2(1+2H)+o(\frac{1}{T^{\gamma}})\notag\\
  &=\int_0^T e^{-2x}x^{2\beta+1}(x+2+\beta) \dif x +(\beta+2)(\beta+1)\int_{0\le x\le z\le T}  e^{-z -x}  x^{\beta}  z ^{\beta} \dif x \dif z \notag\\
  &-\frac12\Gamma^2(\beta)+o(\frac{1}{T^{\gamma}})\notag\\
    &=(2\beta+3)\Gamma(2\beta+2)2^{-(2+2\beta))}+\frac{1+\beta}{2}\Gamma^2(1+\beta)+o(\frac{1}{T^{\gamma}}). 
\end{align}
Moreover, 
the symmetry yields for any $\gamma>0$
\begin{align}
  &  \int_{[0, T]^2 } e^{-z -x}  x^{\beta}  z ^{\beta}\left[ (x\wedge z) -\frac12\right]\dif x \dif z+\int_{0\le x \le z \le T }   e^{-z -x}  x^{\beta}  z ^{\beta} (z-x) \dif x \dif z\notag\\
  &= \int_{0\le x \le z \le T }   e^{-z -x}  x^{\beta}  z ^{\beta} (z+x-1) \dif x \dif z =\int_{[0, T]^2  }   e^{-z -x}  x^{\beta}  z ^{\beta} ( x-\frac12) \dif x \dif z\notag\\
 &=(\beta+\frac12)\Gamma^2(1+\beta)+o(\frac{1}{T^{\gamma}}). \label{dizhongjian}
\end{align}

Plugging \eqref{diyibufen}-\eqref{dizhongjian} into \eqref{j1 biaodashi}, we have
{ \begin{align}
   J_1(T)&= 
   \Gamma^2(1+\beta)\left[\big(\beta+\frac12\big)T - \frac32(1+\beta)^2 +\frac14\right]\notag\\
   &{+}\frac14\left[T^{2\beta}-2\beta T^{2\beta-1}
   +  {\beta}(3\beta-2)T^{2\beta-2} +O(T^{2\beta-3})\right].
   \label{j1 jianjing zhank}
\end{align}}
{ The integration by parts formula implies that 
\begin{align*}
    &\int_{0\le x \le z \le T } e^{x-z} x^{\beta}  z ^{\beta} (x-z+1) \dif x \dif z\\
    &=e^{-T}T^{1+\beta}\int_{0}^T e^x x^{\beta}\dif x-(2\beta+1) \int_{0\le x \le z \le T } e^{x-z} x^{\beta}  z ^{\beta} \dif x \dif z,\\
    &\int_{0\le x \le z \le T } e^{x-z} x^{\beta}  z ^{\beta} \left[z(z-x-1)-\frac12 \right]\dif x \dif z\\
    &=\int_0^T e^x x^{\beta}\dif x \int_x^T z ^{1+\beta} (z-x-1)\dif(-e^{-z}) -\frac12\int_{0\le x \le z \le T } e^{x-z} x^{\beta}  z ^{\beta}  \dif x \dif z\\
    &=(2+\beta) \int_0^{ T } e^{x} x^{\beta}\dif x\int_x^T  z ^{\beta+1}  \dif (-e^{-z})-(1+\beta)\int_0^T e^{-z}z^{\beta}\dif z\int_0^z x^{1+\beta} \dif (e^x)\\
    &-e^{-T}T^{1+\beta}\int_{0}^T e^x x^{\beta}(T-x-1)\dif x - \int_0^T x^{2\beta+1}\dif x-(\beta+\frac32)\int_{0\le x \le z \le T } e^{x-z} x^{\beta}  z ^{\beta}  \dif x \dif z\\
    &=-e^{-T}T^{1+\beta}\int_{0}^T e^x x^{\beta}(T-x+1+\beta)\dif x +(2\beta^2+4\beta+\frac32)\int_{0\le x \le z \le T } e^{x-z} x^{\beta}  z ^{\beta}  \dif x \dif z.
\end{align*}
The above two identities together with the symmetry imply that }
\begin{align}
    J_2(T)&=\int_{0\le x \le z \le T } e^{x-z} x^{\beta}  z ^{\beta}\left[-(z-x)(T-z)+T-z-\frac12 \right]  \dif x \dif z\notag\\
    &=\left(-(2\beta+1)T+2(1+\beta)^2-\frac12 \right)\int_{0\le x \le z \le T } e^{x-z} x^{\beta}  z ^{\beta} \dif x \dif z\notag\\
    &+e^{-T}T^{1+\beta}\int_{0}^T e^x x^{\beta}(x-1-\beta)\dif x:=J_{21}(T)+J_{22}(T).\label{j2t kaishi}
\end{align}
Since $\beta=2H-1$, it is clear that 
 $J_{21}(T)=0$ when $H=\frac14$. 
It follows from Lemma~\ref{asymptotic expansion key} that when $H\in(0,\frac14)\cup(\frac14, \frac34)$,
\begin{align}
  J_{21}(T)&= 
 \frac{(1-4H)\Gamma^2(2H)}{2\cos(2H\pi)} [T-\frac{4H+1}{2}] -T^{4H} +4HT^{4H-1}+O(T^{4H-4})\notag\\
 &+(1-4H)\left[H+\frac14 +\frac{\beta(\beta-1)}{2\beta-1}\right]T^{4H-2}+{\frac{\beta(2\beta+1)(4\beta^2-4\beta-1)}{2(2\beta-1)}}T^{4H-3}. \label{j21t}
\end{align}
Using the integration by parts and Lemma~\ref{asymptotic expansion key-0}, we have 
\begin{align}
    J_{22}(T)&=e^{-T}T^{1+\beta}\left[\int_0^T x^{1+\beta}\dif e^x -(1+\beta)\int_0^T e^x x^{\beta}\dif x\right]\notag\\
    &=T^{4H}-4H \left[T^{4H-1}-(2H-1) T^{4H-2}+(2H-1)(2H-2)T^{4H-3} +O(T^{4H-4})\right].\label{j22t}
\end{align}
Substituting \eqref{j21t} and \eqref{j22t} into \eqref{j2t kaishi}, we have
\begin{align}
    J_2(T)
&= \frac{(1-4H)\Gamma^2(2H)}{2\cos(2H\pi)}\Big[T-\frac{4H+1}{2} \Big] \notag\\
&+\left[4H(2H-1)+(1-4H)\left[H+\frac14 +\frac{\beta(\beta-1)}{2\beta-1}\right]\right]T^{4H-2}+O(T^{4H-3})\notag\\
&=\frac{(1-4H)\Gamma^2(2H)}{2\cos(2H\pi)}\Big[T-\frac{4H+1}{2} \Big]+\left[-\frac14+\frac{1}{2(4H-3)}\right] T^{4H-2}\notag\\
&+\left[H-\frac{3}{2}-\frac{1}{4H-3}\right]T^{4H-3} +O(T^{4H-4}).\label{j2T biaods}
    \end{align}
Substituting the computation results \eqref{j1 jianjing zhank} and \eqref{j2T biaods} into  \eqref{i1Tfenjie}, we obtain \eqref{fenjie 001}.
\end{proof}

\begin{lemma}\label{L2T calculat}
Suppose $\beta=2H-2$ and 
\begin{align*}
    L_2(T)&:=\int_{0}^T e^{-t_1}\dif t_1 \int_{[0,t_1]^3}  {e^{ s_1-|t_2-s_2|}} (t_1-t_2)^{\beta}   \abs{s_1-s_2}^{\beta}  \dif t_2 \dif \vec{s}.
\end{align*}
 When $H\in (\frac12,\frac34)\cup(\frac34,1)$, the following expansion holds:
\begin{align}
  L_2(T)
  &=\frac12\Gamma^2({1+\beta})\left[(2\beta+3)(1-\frac{1}{\cos(\beta \pi)})(T-\frac{2\beta+3}{2})-(\beta+1)(\beta+3) \right]\notag\\
    &+\frac{1}{(\beta+1)(2\beta+1)}T^{2\beta+2}- \frac{2}{2\beta+1}T^{2\beta+1}+O(T^{2\beta}). \label{fenjie 002-l2t-1}
\end{align}
When $H=\frac34$, the following expansion holds:
\begin{align}
  L_2(T)=2\left[(T-1)\log T+(c+\frac{\pi-3}{2})T+1-c-\frac{13\pi}{16}+\frac{1}{T}\right]+O(T^{-2}). \label{34teshuqkzk}  
\end{align}
\end{lemma}
\begin{proof}
The domain $[0,T]\times [0,t_1]^3$ can be divided into six disjoint regions according
to the distinct orders of $t_2,s_1,s_2$ as \eqref{domain divided},
 and for $i=1,2,3,1',2',3'$, we denote
$$L_{2i}(T)=\int_{\Delta_i}  {e^{ s_1-t_1-|t_2-s_2|}} (t_1-t_2)^{\beta}   \abs{s_1-s_2}^{\beta}  \dif \vec{t} \dif \vec{s}.$$   
So $L_2(T) = \sum_{i}L_{2i}(T)$. Making substitution $x=\abs{s_1-s_2}, y=\abs{(s_1\vee s_2)-t_2}, z=t_1-t_2$, similar to \eqref{diyibufen-01}, \eqref{dierbufen}, and \eqref{disanbufen}, we have that
\begin{align}
  L_{21}(T)+ L_{21'}(T)  &= \int_{0\le x\le y\le z\le t_1\le T }  \left({e^{-z -x}+e^{x-z}}\right) x^{\beta}  z ^{\beta} \dif x \dif y \dif z\dif  t_1\notag\\
   &=\int_{0\le x \le z \le T }  \left({e^{-z -x}+e^{x-z}}\right) x^{\beta}  z ^{\beta} (z-x)(T-z) \dif x \dif z,\label{L2diyibufen-01}
\end{align}
and 
\begin{align}
  L_{22}(T)+ L_{22'}(T)  &= \int_{0\le y\le x,z\le x+z-y\le t_1\le T }  \left({e^{-z -x}+e^{-z-x+2y}}\right) x^{\beta}  z ^{\beta} \dif x \dif y \dif z\dif  t_1,\label{L2dierbufen-01}
\end{align}
and 
\begin{align}
  L_{23}(T)+ L_{23'}(T)  &= 2\int_{[0,T]^4,\,0\le x+y+z\le  \le t_1\le T }   e^{-z-x-2y}  x^{\beta}  z ^{\beta} \dif x \dif y \dif z\dif  t_1. \label{L2dierbufen-01}
\end{align}

Using the expression \eqref{dierbufen}, we have that
\begin{align}
  L_{21}(T)+ L_{21'}(T)+   L_{22}(T)+ L_{22'}(T) :=J_1(T)+\bar{J}_2(T), \label{j1jwbarT}
\end{align}where $J_1(T)$ is given by \eqref{j1 biaodashi}, and $\bar{J}_2(T)$ is given as follows:
\begin{align}
    \bar{J}_2(T)&:=\int_{0\le x \le z \le T } e^{x-z} x^{\beta}  z ^{\beta}\left[ (z-x)(T-z) +T- z -\frac12\right]\dif x \dif z\notag\\
    &=T\left[\int_{0\le x \le z \le T } e^{x-z} x^{\beta}  z ^{\beta}  (z-x+1) \dif x \dif z \right]-\int_{0\le x \le z \le T } e^{x-z} x^{\beta}  z ^{\beta+1}  (z-x+1) \dif x \dif z\notag\\
    &-\frac12\int_{0\le x \le z \le T } e^{x-z} x^{\beta}  z ^{\beta} \dif x \dif z:=T\bar{J}_{21}(T)+\bar{J}_{22}(T)-\frac12\bar{J}_{23}(T),\label{j2 fenjieshizi}
\end{align}where $\bar{J}_{23}(T)=\int_{0\le x \le z \le T } e^{x-z} x^{\beta}  z ^{\beta} \dif x \dif z$. The integration by parts implies that 
\begin{align}
    \bar{J}_{21}(T)&=\int_0^T e^{x}x^{\beta}\dif x \int_x^T z^{1+\beta}\dif(-e^{-z})-\int_0^T e^{-z}z^{\beta}\dif z \int_0^z x^{1+\beta}\dif e^{x}+\bar{J}_{23}(T)\notag\\
    &=-T^{1+\beta}e^{-T}\int_0^T e^x x^{\beta}\dif x +(2\beta+3)\bar{J}_{23}(T),\label{j21t zhankai}
\end{align}and 
\begin{align}
    \bar{J}_{22}(T)&=\int_0^T e^{x}x^{\beta}\dif x \int_x^T z^{1+\beta}(z-x+1)\dif(e^{-z})\notag\\
    &=\int_0^T e^{x}x^{\beta}\dif x \left[ e^{-z}z^{1+\beta}(z-x+1)|_{z=x}^T -\int_x^T e^{-z}\left[(2+\beta)z^{\beta+1} +(1-x)(1+\beta)z^{\beta}\right]\dif z\right]\notag\\
    &=T^{2+\beta}e^{-T}\int_0^T e^x x^{\beta}\dif x-\frac{\beta+2}{\beta+1}T^{2\beta+2}+2(2+\beta)T^{1+\beta}e^{-T}\int_0^T e^x x^{\beta}\dif x\notag\\
    &-2(1+\beta)(2+\beta)\bar{J}_{23}(T).\label{j22t zhankai}
\end{align} Substituting \eqref{j21t zhankai} and \eqref{j22t zhankai} into \eqref{j2 fenjieshizi} and using Lemmas~\ref{asymptotic expansion key-0} and~\ref{asymptotic expansion key}, we have that when $\beta\neq -\frac{1}{2}$,
\begin{align}
  \bar{J}_{2}(T)&=\left[(2\beta+3)T-2(1+\beta)(2+\beta)-\frac12\right] \bar{J}_{23}(T)\notag\\
  &-\frac{\beta+2}{\beta+1}T^{2\beta+2}+ 2(2+\beta) T^{1+\beta}\left[T^{\beta}-\beta T^{\beta-1}+O(T^{\beta-2})\right]\notag\\
  &=\frac{1}{(\beta+1)(2\beta+1)}T^{2\beta+2}- \frac{2}{2\beta+1}T^{2\beta+1}+O(T^{2\beta})\notag\\
  &-\frac{\Gamma^2(1+\beta)}{2\cos(\pi \beta)}\left[(2\beta+3)T-2(1+\beta)(2+\beta)-\frac12\right].\label{jbar2Thuajian}
\end{align} Combining \eqref{j1jwbarT}, \eqref{jbar2Thuajian} with \eqref{j1 jianjing zhank}, we have 
\begin{align}
    & L_{21}(T)+ L_{21'}(T)+   L_{22}(T)+ L_{22'}(T) \notag\\
    &=\frac{1}{(\beta+1)(2\beta+1)}T^{2\beta+2}- \frac{2}{2\beta+1}T^{2\beta+1} -\frac{\Gamma^2(1+\beta)}{2\cos(\pi \beta)}(2\beta+3)\left[T-\frac{2\beta+3}{2}\right]\notag \\
  &+ \frac12\Gamma^2(1+\beta)\left[\big(2\beta+1\big)T - 3(1+\beta)^2+\frac12\right]+O(T^{2\beta}).\label{l21_24}
\end{align}

Finally, we have 
\begin{align*}
   &L_{23}(T)+ L_{23'}(T)\notag\\
   &=\int_{[0,T]^2,x+z\le T}   e^{-z-x}  x^{\beta}  z ^{\beta} \dif x \dif z \int_0^{T-x-z}\dif (- e^{-2y}) \int_{0}^{T-x-y-z} \dif  t_1\notag\\
   &=\int_{[0,T]^2,x+z\le T}   e^{-z-x}  x^{\beta}  z ^{\beta} \dif x \dif z  \left[-(T-x-y-z)e^{-2y}|_{y=0}^{T-x-z} -\int_{0}^{T-x-z} e^{-2y}\dif y \right]\notag\\
   &=\int_{[0,T]^2,x+z\le T}   e^{-z-x}  x^{\beta}  z ^{\beta}  \left[ T-x-z-\frac12 +\frac12 e^{-2(T-x-z)}\right] \dif x \dif z.
\end{align*}
Since for any $\gamma>0$,
\begin{align*}
    0<\int_{[0,T]^2,x+z\le T}   e^{-z-x-2(T-x-z)}  x^{\beta}  z ^{\beta}   \dif x \dif z<\frac{1}{(1+\beta)^2}e^{-T}T^{2(1+\beta)}=o(\frac{1}{T^{\gamma}}),
\end{align*}and 
\begin{align*}
    0&<-\int_{[0,T]^2,x+z> T}   e^{-z-x}  x^{\beta}  z ^{\beta}  \left[ T-x-z-\frac12  \right] \dif x \dif z\\
    &< e^{-T}\int_{[0,T]^2,x+z> T}   x^{\beta}  z ^{\beta}  \left[ x+z+\frac12  \right] \dif x \dif z=o(\frac{1}{T^{\gamma}}),
\end{align*} we have that for any $\gamma>0$,
\begin{align}
    L_{23}(T)+ L_{23'}(T)&=\int_{[0,T]^2 }   e^{-z-x}  x^{\beta}  z ^{\beta}  \left[ T-x-z-\frac12  \right] \dif x \dif z +o(\frac{1}{T^{\gamma}})\notag\\
    &=\Gamma^2({1+\beta})(T-\frac12)-2\Gamma({1+\beta})\Gamma({2+\beta})+o(\frac{1}{T^{\gamma}})\notag\\
    &=\frac12\Gamma^2({1+\beta})(2T-4\beta-5)+o(\frac{1}{T^{\gamma}}). \label{l3l3prime}
\end{align}

Combining \eqref{l3l3prime} with \eqref{l21_24}, we have 
\begin{align*}
    L_2(T) &= \sum_{i=1}^3 [L_{2i}(T)+L_{2i'}(T)] \\
    &=\frac12\Gamma^2({1+\beta})\left[(2\beta+3)(1-\frac{1}{\cos(\beta \pi)})(T-\frac{2\beta+3}{2})-(\beta+1)(\beta+3) \right] \\
    &+\frac{1}{(\beta+1)(2\beta+1)}T^{2\beta+2}- \frac{2}{2\beta+1}T^{2\beta+1}+O(T^{2\beta}).
\end{align*}

Finally, \eqref{34teshuqkzk} can be obtained from Lemma~\ref{asymptotic expansion key} in the same vein.
\end{proof}


\vspace{1cm} \noindent {\bf {References}} \small

\begin{itemize}

\item [{[1]}] F Alazemi, A Alsenafi, Y Chen, H Zhou.  \emph{Parameter Estimation for the Complex Fractional Ornstein-Uhlenbeck Processes with Hurst parameter $H\in (0,\,\frac12)$},   { Chaos Solitons Fractals}, 2024, 188: 1155562024.

\item [{[2]}]M Balde, R Belfadli, K Es-Sebaiy.  \emph{Kolmogorov bounds in the CLT of the LSE for Gaussian
Ornstein-Uhlenbeck processes},  { Stoch. Dynam.} 2023, 23(4) 2350029.

\item [{[3]}]  Y Chen, X Gu. \emph{ An Improved Berry-Esséen Bound of Least Squares Estimation for Fractional Ornstein-Uhlenbeck Processes}, { Acta Mathematica Scientia, Series A}, 2023, 43(3): 855-882. (in Chinese)

\item [{[4]}]  Y Chen, Y Li.  {\it Berry-Ess\'een bound for the parameter estimation of fractional Ornstein-Uhlenbeck processes with the hurst parameter $H \in (0, \frac12)$},  Commun. Stat. Theory Methods, 2021, {  50}(13), 2996-3013.

\item [{[5]}] Y Chen, Y Li, H Zhou.  {\it Berry-Ess\'en bounds for the statistical estimators of an Ornstein-Uhlenbeck process driven by a general Gaussian noise},  Fract. Calc. Appl. Anal, 2025, 28: 2607-2637.

\item [{[6]}] Y Chen, H Zhou.  {\it Parameter estimation for an Ornstein-Uhlenbeck process driven by a general gaussian noise},  Acta Math. Sci. Ser. B (Engl. Ed.), 2021, 41(2) 573-595.

\item [{[7]}]P Cheridito, H Kawaguchi, M Maejima.  {\it Fractional Ornstein-Uhlenbeck processes},  Electron. J. Probab. 2003, 8: 1-14.

\item [{[8]}]  S Douissi, K Es-Sebaiy,  G Kerchev, I Nourdin. \emph{ Berry-Ess\'een bounds of second moment estimators for Gaussian processes observed at high frequency}, Electron J Stat,  2022, {16}(1), 636-670

\item [{[9]}]K Es-Sebaiy, F Alazemi.  \emph{New Kolmogorov bounds in the CLT for random ratios and applications}, {Chaos Solitons Fractals}, 2024, 181: 114686.

\item [{[10]}] Y Hu, D Nualart.  {\it Parameter estimation for fractional Ornstein-Uhlenbeck processes},  Statist Probab Lett, 2010, {  80}(11-12).

\item [{[11]}] Y Hu, D Nualart, H Zhou.  {\it Parameter estimation for fractional Ornstein-Uhlenbeck
processes of general Hurst parameter}, Statistical Inference for Stochastic Processes, 2019, 22(1): 111–42. 

\item [{[12]}] D Nualart. The Malliavin calculus and related topics. Springer.  2006.

\end{itemize}\vskip 10mm

\noindent School of Big Data, Baoshan University, Baoshan, 678000, Yunnan, China.\\
\indent  Email: zhishi@pku.org.cn

\end{document}